\documentclass[10pt,twocolumn]{article}
\pdfoutput=1
\usepackage{microtype}
\usepackage[pdftex,colorlinks,linkcolor=black,urlcolor=blue]{hyperref}
\usepackage[flushleft]{threeparttable}

\newcommand{\footlabel}[2]{%
    \addtocounter{footnote}{1}%
    \footnotetext[\thefootnote]{%
        \addtocounter{footnote}{-1}%
        \refstepcounter{footnote}\label{#1}%
        #2%
    }%
    $^{\ref{#1}}$%
}
\newcommand{\footref}[1]{%
    $^{\ref{#1}}$%
}

\newcommand\blfootnote[1]{%
  \begingroup
  \renewcommand\thefootnote{}\footnote{#1}%
  \addtocounter{footnote}{-1}%
  \endgroup
}

\newcommand{\subs}[1]{\smallskip\noindent\textbf{#1}}

\begin{document}

\twocolumn[%
\centerline{\Large\bf Mathematicians Take a Stand}
\centerline{Douglas N.\ Arnold and Henry Cohn}
\medskip]

Mathematicians care deeply about the mathematical literature.
We devote much of our lives to learning from it, expanding it,
and guaranteeing its quality.  We depend on it for our
livelihoods, and our contributions to it will be our
intellectual legacy.\blfootnote{$\!\!\!\!$Douglas N.\ Arnold is
McKnight Presidential Professor of Mathematics at the
University of Minnesota. His email address is
\href{mailto:arnold@umn.edu}{arnold@umn.edu}.

Henry Cohn is
principal researcher at Microsoft Research New England and
adjunct professor of mathematics at MIT.  His e-mail address is
\href{mailto:cohn@microsoft.com}{cohn@microsoft.com}.}

It has long been anticipated that technological advances will
make the literature more affordable and accessible. Sadly, this
potential is not being fully realized. The prices libraries
pay for journals have been growing with no end in sight,
even as the costs of publication and distribution have gone
down, and many libraries are unable to maintain
their subscriptions.%
\footnote{For example, MIT's spending on serials
increased by 426 percent over the period 1986--2009, while
the number of serials purchased decreased by 16 percent,
and the Consumer Price Index increased by only 96 percent.}

The normal market mechanisms we count on to keep prices in
check have failed for a variety of reasons.  For example,
mathematicians have a professional obligation to follow the
relevant literature, which leads to inelastic pricing. This
situation is particularly perverse because we provide the
content, editorial services, and peer review free of charge,
implicitly subsidized by our institutions.  The journal
publishers then turn to the same institutions and demand prices
that seem unjustifiable.

Although the detailed situation is complex, the fundamental cause of
this sad state of affairs is not hard to find. While libraries are
being forced to cut acquisitions, a small number of commercial
publishers have been making breathtaking profits year after year. The
largest of these, Elsevier, made an adjusted operating profit of \$1.12 billion in 2010 on
\$3.14 billion in revenue, for a profit margin of 36 percent, up from 35 percent in
2009
and 33 percent in 2008.%
\footnote{Reed Elsevier Annual Report 2010, SEC form 20-F (based on data
from p.~25 and average exchange rate from p.~6).}
Adding insult to injury, Elsevier has aggressively pushed
bundling arrangements that result in libraries paying for
journals they do not want and that obscure the actual
costs.\footnote{T.~Bergstrom, Librarians and the terrible
fix: economics of the Big Deal, \emph{Serials} \textbf{23} (2010),
77--82.} They have fought transparency of pricing, going so far
as to seek a court injunction in an unsuccessful attempt to
stop a state university from revealing the terms of their
subscription contract.  They have imposed
unacceptable restrictions on dissemination by authors. And,
while their best journals make important contributions to the
mathematical literature, Elsevier also publishes many weaker
journals, some of which have been caught in major lapses of
peer review or ethical standards.  These scandals have done
harm to the integrity and reputation of mathematics.

This situation has been extensively analyzed many times before,
including in the \emph{Notices}.  There have been some high-profile actions,
such as mass resignations of entire Elsevier editorial boards over
pricing concerns: the \emph{Journal of Logic Programming} in 1999, the \emph{Journal
of Algorithms} in 2003, and \emph{Topology} in 2006.  These boards have done a
valuable service for the community by founding replacement journals,
but there has been little relief from the overall trend. As Timothy
Gowers wrote in his blog in January, ``It might seem inexplicable that
this situation has been allowed to continue. After all, mathematicians
(and other scientists) have been complaining about it for a long time.
Why can't we just tell Elsevier that we no longer wish to publish with
them?'' Gowers then revealed that he had been quietly boycotting
Elsevier for years, and he suggested it would be valuable to create a
website where like-minded researchers could \emph{publicly} declare
their unwillingness to contribute to Elsevier journals.

Within days, Tyler Neylon responded to this need by creating
\url{http://thecostofknowledge.com}.  More than 2,000 people
signed on in the first week, and participation has
grown steadily since then, to over 8,000 as of early March. Each
participant chooses whether to refrain from publishing papers
in, refereeing for, or editing Elsevier journals.  The boycott
is ongoing, and it holds the promise of sparking real change.
\emph{We urge you to consider adding your voice.}

The boycott is a true grassroots movement.  No individual or
group is in charge, beyond Gowers's symbolic position
as the first boycotter.  However, a group of thirty-four mathematicians%
\footnote{Scott Aaronson, Douglas N.\ Arnold, Artur Avila, John
Baez, Folkmar Bornemann, Danny Calegari, Henry Cohn, Ingrid
Daubechies, Jordan Ellenberg, Matthew Emerton, Marie Farge,
David Gabai, Timothy Gowers, Ben Green, Martin Gr\"otschel,
Michael Harris, Fr\'ed\'eric H\'elein, Rob Kirby, Vincent
Lafforgue, Gregory F.\ Lawler, Randall J.\ LeVeque, L\'aszl\'o
Lov\'asz, Peter J.\ Olver, Olof Sisask, Terence Tao, Richard
Taylor, Bernard Teissier, Burt Totaro, Lloyd N.\ Trefethen,
Takashi Tsuboi, Marie-France Vign\'eras, Wendelin Werner, Amie
Wilkinson, and G\"unter M.\ Ziegler.} (including Gowers and the
authors of the present article) issued their best attempt at a
consensus statement of purpose for the boycott.  It is
available online,%
\footnote{See \url{http://umn.edu/~arnold/sop.pdf} or the
March 2012 \emph{London Mathematical Society Newsletter}.} and we
highly recommend it for reading.  For reasons of space, we will not
cover every aspect of that statement here.

Before we proceed, we must address two pressing questions about
the boycott. First, why is a boycott appropriate? After all,
Elsevier employs many reasonable and thoughtful people, and
many mathematicians volunteer their services, helping to
produce journals of real value. Isn't a boycott overly
confrontational?  Could one not take a more collaborative
approach? Unfortunately, such approaches have been tried time
and again without success.  Fifteen years of reasoned
discussions have failed to sway Elsevier.%
\footnote{R.~Kirby, Comparative prices of math journals,
1997, \url{http://math.berkeley.edu/~kirby/journals.html};
J.~Birman, Scientific publishing: a mathematician's
viewpoint, \emph{Notices of the AMS} \textbf{47} (2000), 770--774;
R.~Kirby, Fleeced?, \emph{Notices of the AMS} \textbf{51}
    (2004), 181;
W.~Neumann, What we can do about journal pricing, 2005,
    \url{http://www.math.columbia.edu/~neumann/journal.html};
D.~N.~Arnold, Integrity under attack: the state of
scholarly publishing, \emph{SIAM News} \textbf{42} (2009), 2--3;
P.~Olver, Journals in flux, \emph{Notices of
    the AMS} \textbf{58} (2011),
    1124--1126.}
Elsevier's leadership seems to be driven only by their fiduciary
responsibility to maximize profit for their shareholders.  The
one hope we see for change is to demonstrate that their
business depends on us and that we will not cooperate with
them unless they earn our respect and goodwill.

Second, why is the focus solely on Elsevier?  Some of the problems we
discuss are common among large commercial publishers, and indeed we hope the
boycott will help spur changes in the whole industry.  But we must
start somewhere, and we believe it is more effective to focus on one
publisher whose behavior has been particularly egregious than to
directly confront an entire industry at once.  Many of the successful
boycotts in history took the same tack.

\begin{table*}[ht]
\begin{center}
\begin{threeparttable}
\caption{Summary information for six journals.}
\label{table:summary}
\begin{tabular}{lllrrrr}
\hline
\hline\vrule height4mm width0mm
Journal & Publisher & Metrics         & Price & \$/art.& \$/page & \$/cite \\[1mm]
\hline\vrule height6mm width0mm
Annals of Mathematics & Princeton      & 3.7/A* &   \$447 &  5.39  & 0.12 & 0.06 \\
SIAM J.~Appl.~Math.   & SIAM     & 1.8/A* &   \$642 &  5.95 & 0.27 & 0.13 \\
Journal of the AMS    & AMS      & 3.6/A* &   \$300 &  9.09 & 0.24 & 0.13 \\
Advances in Mathematics & Elsevier & 1.6/A* & \$3,899 & 11.53 & 0.35 & 0.90 \\
Journal of Algebra      & Elsevier & 0.7/A* & \$6,944 & 13.89 & 0.75 & 1.22 \\
Journal of Number Theory  & Elsevier & 0.6/B & \$2,745 & 17.49 & 1.12 & 1.91 \\[2mm]
\hline\hline
\end{tabular}
\begin{tablenotes}
\footnotesize\item[\quad]
    Metrics are the 2010 5-year impact factor from
    \emph{Journal Citation Reports} and the 2010 rating by the
    Australian Research Council (based on expert opinion).
    A* = top-rated, B = ``solid, though not outstanding''.

\item[\quad] Elsevier prices are the amounts
    actually paid by the University of Minnesota for electronic-only institutional
    subscriptions in 2012.  The lowest prices we could find on the Elsevier website as of
    February 29 were \$3,555.20, \$5,203, and \$2,226.40.
    The \emph{Annals} price is again the actual amount paid by UMN, which
    is slightly greater than
    the \$435 list price.  The SIAM and AMS prices are the list prices, although UMN paid
    less because of institutional membership.

\item[\quad] Columns 5--7 normalize by the most recent data
    available: the numbers of articles and pages published in 2011 and the number of citations to the journal
    made in 2010 (as reported in \emph{Journal Citation Reports}).\vspace{-2mm}
\end{tablenotes}
\end{threeparttable}
\end{center}
\end{table*}

\subs{Journal Pricing.} Table~\ref{table:summary} exhibits prices for
three of Elsevier's mathematics journals: \emph{Advances in Mathematics}, the
\emph{Journal of Algebra}, and the \emph{Journal of Number Theory}. For comparison,
the table includes three more affordable journals.

The \emph{Annals of Mathematics}, published by the Princeton math department and IAS,
provides exceptional
quality at a rock-bottom price that just covers costs. The other
two are highly regarded journals published
by the Society for Industrial and Applied Mathematics (SIAM) and
by the American Mathematical Society (AMS). Both of these organizations make a
profit on their journal publishing operations, which helps to
subsidize their other activities. For example, in 2011 SIAM's
journal publication costs, including overhead, were 89 percent
of their subscription revenues, resulting in an 11 percent profit
margin.

Elsevier's recent pricing changes, apparently in response to
the boycott, have at times led to multiple conflicting prices
on their website. We have listed the prices actually paid by
the University of Minnesota in 2012, but the notes after the
table indicate the lowest prices we found offered on the Web.
We made no attempt to select the highest-priced Elsevier
journals, and in fact \emph{Advances in Mathematics} is among the most
affordable.  For comparison, the 2011 prices per page of the
thirty-six Elsevier journals listed in the AMS journal price
survey ranged from \$0.33 (\emph{Advances in Mathematics}) to \$4.05
(\emph{Mathematical Social Sciences}), with a mean of \$1.35 and a
median of \$0.96.

As shown in the table, the prices of the SIAM and AMS journals are
within a factor of two of that of the \emph{Annals}, with differences
depending on whether one normalizes the raw journal price by number of
articles, pages, or citations.  But the Elsevier prices are a different
story. The price per page of the \emph{Journal of Algebra}, for example, is
triple that of the society journals and six times that of the \emph{Annals},
and the \emph{Journal of Number Theory} is 50 percent more expensive yet.

\begin{table*}
\begin{center}
\begin{threeparttable}
\caption{Historical prices per page in constant 2012 dollars.}
\label{table:historical}
\begin{tabular}{lrrrrrrrr}
\hline
\hline\vrule height4mm width0mm
 & 1994 & 1997 & 2000 & 2003 & 2006 & 2009 & 2010 & 2011 \\[1mm]
\hline\vrule height6mm width0mm
Annals of Mathematics\hspace{1cm} & 0.19 & 0.20 & 0.14 & 0.15 & 0.13 & 0.13 & 0.09 & 0.10\\
SIAM J.~Appl.~Math. & 0.20 & 0.24 & 0.23 & 0.25 & 0.27 & 0.24 & 0.18 & 0.27\\
Journal of the AMS & 0.22 & 0.26 & 0.27 & 0.29 & 0.30 & 0.27 & 0.25 & 0.24\\
Advances in Mathematics & 0.65 & 0.74 & 0.95 & 1.01 & 0.55 & 0.61 & 0.44 & 0.33\\
Journal of Algebra & 0.36 & 0.43 & 0.50 & 0.73 & 0.60 & 0.77 & 0.92 & 0.66\\
Journal of Number Theory & 0.57 & 0.67 & 0.98 & 1.01 & 1.04 & 0.86 & 0.95 & 1.05\\[2mm]
\hline\hline
\end{tabular}
\begin{tablenotes}
\footnotesize\item[\quad] Prices are from the AMS journal price
survey
(\url{http://www.ams.org/membership/mem-journal-survey}),
adjusted using the Consumer Price Index. \vspace{-5mm}
\end{tablenotes}
\end{threeparttable}
\end{center}
\end{table*}

Moreover, as demonstrated in Table~\ref{table:historical}, this
problem has grown over time. The inflation-adjusted prices per
page of the \emph{Journals of Algebra} and \emph{Number Theory} increased by
more than 80 percent between 1994 and 2011, compared with much
smaller increases for the society journals and a decrease for
the \emph{Annals}. It is noteworthy that the recent prices of \emph{Advances
in Mathematics}, while still high, have come closer to the
prices of the society journals. This supports our belief that
Elsevier could offer substantially lower prices and still make
a reasonable profit.

We do not mean to suggest that publishing is cheap in the
electronic age. True, electronic distribution is very cheap:
the arXiv requires
just \$7 per submission, or 1.4 cents per download, in funding.%
\footnote{\url{http://arXiv.org/help/support/faq}.} But journal
publishing involves significant additional costs, such as IT
infrastructure, administrative support, oversight, sales,
copyediting, typesetting, archiving, etc.  Many of these costs
scale roughly with the number of published pages, and some of
them benefit from economies of scale (so large publishers like
Elsevier should, if anything, achieve lower costs).

Of course, journals are not all the same. A low-circulation journal may
need to command a higher price per page to stay afloat.  The community
might find some such journals too expensive to support, but one viewed
as worthy of support might reasonably charge a higher price until more
libraries subscribe.  Another journal might have extraordinary
expenses, for example from translation.  But these factors do not apply
to the cases we have considered or to many other Elsevier
journals.

We see no good reason to pay much more for Elsevier journals
than for journals that earn mathematical societies a tidy profit. Even the
price targets for mathematics journals that Elsevier
announced in response to the boycott%
\footlabel{fn}{D.~Clark and L.~Hassink, A letter to
    the mathematics community, February 27, 2012,
    \url{http://www.elsevier.com/wps/find/P11.cws_home/lettertothecommunity}.
}---\$0.50 to \$0.60 per page---would leave their journals costing twice as much as the
comparison journals in Table~\ref{table:summary}. Elsevier's prices
have become far out of proportion and have a way to go to return to
reasonable.

\subs{What's the Big Deal about Bundling?} Bundling refers to
grouping together collections of journals and selling access as
a single product, discounted from list price.   Elsevier commonly
negotiates bundles including all the journals to which the library has recently subscribed.
The bundles may also include access to nearly all of
Elsevier's roughly 2,000 journals. Librarians have termed
enormous bundles ``the big deal''.

While there is nothing wrong with offering quantity discounts,%
\footnote{For example, Mathematical Sciences Publishers offers
a bundle of six mathematics journals at a 31 percent discount,
bringing their price down to \$0.08 per page.} it is the way in
which Elsevier and other large publishers have implemented
bundling that is objectionable. They have turned it into a
powerful tool for subverting the market forces that would keep
prices in check. The then director of Harvard's library
summarized it thus: ``Elsevier is among a handful of journal
publishers whose commercial bundling practices are squeezing
library budgets. Their licensing programs require libraries to
maintain large, fixed levels of expenditure, without the
ability to cancel unneeded
subscriptions.''%
\footnote{\url{http://hul.harvard.edu/news/2004_0101.html}.}

Let us see how this works. While Elsevier has gone to great
lengths to keep the details of their bundle contracts secret, some have
come to light, thanks to open records laws.%
\footnote{T.~Bergstrom, P.~Courant, and R.~P.~McAfee, \emph{Big Deal
Contract Project}.}  Judging by the contracts we have seen and
librarians we have consulted, it works essentially as follows.

The university commits to subscribing to the
journals it currently receives at a negotiated total price
that is typically around the same as they were previously paying and to
continuing to subscribe to them for a period of three to five years with
annual price increases. Elsevier has called this the ``Complete
Collection'', and it is a large expense. For example, for the
University of Minnesota in 2006 it came to \$1.8 million (about
18 percent of their total serials budget), and for the University of
Michigan in 2007 it was \$1.9 million.  In both cases, 5 percent
yearly price increases were built into the contracts, although
the actual rate of inflation for the contract periods was only
about 2 percent per year. Cancellation of titles in the Complete
Collection is restricted, which makes it difficult or impossible
to cut back on the expenditure.

For an additional fee Elsevier offers their ``Freedom Collection'',
which adds deeply discounted access to nearly all of the Elsevier
journals to which the library had \emph{not} chosen to subscribe.  This
option cost the University of Michigan about \$19,000 more in 2007, inflating 5 percent a
year thereafter. The University of Minnesota elected against it.

Although prices increase quickly inside the bundle, list prices
can increase even more quickly, so a university that decides
not to renew its bundle may face a steep price increase to hold
onto the journals it wants. Because of bundling, ever
larger portions of library budgets are locked into Elsevier
contracts, budgetary pressures force the
cancellation of titles from smaller publishers,
and funds for new subscriptions disappear.  Furthermore, bundling leads to a lack
of clarity on pricing. The discounts on the additional journals
in the Freedom Collection can sound impressive, but it is the
pricing of the primary subscriptions that drains library
budgets.

The constraints imposed by bundles have led some universities
to conclude that even paying exorbitant prices for the journals
they choose is a better deal.  Harvard, MIT, the University of
Minnesota, and others have now gone this route.  However, many
academic libraries remain tied to the big deal.

Price disclosure is necessary for a well-functioning market
with competitive pricing, so the lack of transparency in
bundling contracts is particularly troubling.  As an Elsevier
vice president wrote in support of Elsevier's 2009 lawsuit
enjoining Washington State University from revealing the prices
of their subscriptions, ``Elsevier representatives apply
pricing formulae and methods which are not generally known (to
our competitors or potential customers)'' and ``disclosure
could disadvantage Elsevier in that, if its pricing to customer
X was known to
customers Y and Z, the latter could demand the same pricing''.%
\footnote{\url{http://www.econ.ucsb.edu/~tedb/Journals/WSUCourtCase/ElsevierStatementbySalesChief.pdf}.}
Elsevier may indeed profit from keeping Y and Z in the dark, but the
academic community values sunlight.  Without transparency of
subscription contracts and costs, the community will remain skeptical
of Elsevier's pricing, whatever changes they make to list prices.

\subs{Posting Policies.} Gowers's suggestion of an Elsevier
boycott struck a chord in many researchers. Besides pricing and
bundling, there are other issues that have contributed to so
many researchers' readiness to abandon Elsevier. One of these
is Elsevier's policies concerning dissemination.   Thanks to the Internet,
authors have additional ways of disseminating their work besides the
printed journal and journal website. For example, it has become common practice
for authors to post a finalized version of their manuscript to
a repository such as the arXiv for open dissemination,
as allowed by many publishers.\footnote{{K.~Fowler, Do
mathematicians get the
    author rights they want?, \emph{Notices of the AMS}
    \textbf{59} (2012), 436--438.}} Elsevier's actions suggest that
they view this development primarily as a threat to their
profits, not as an opportunity to advance mathematics or
increase their authors' readership. In short, their
interests are not aligned with ours.

Elsevier's policies are complex and difficult to understand. In the
words of the scholarly communications officer at Duke University,
``it seems clear that the intent of these statements, policies and
contracts is not to clarify the authors' obligations so much as it is
to confuse and
intimidate them.''%
\footnote{K.~Smith, What a mess!, \emph{Scholarly
    Communications @ Duke}, July 7, 2011,
    \url{http://blogs.library.duke.edu/scholcomm/2011/07/07/what-a-mess/}.}
Their posting policy%
\footnote{\url{http://www.elsevier.com/wps/find/authorsview.authors/postingpolicy},
accessed March 3, 2012.}
 specifically prohibits posting an
``accepted author manuscript''---the author's own version of a
manuscript that has been accepted for publication---on an
e-mail list, a subject repository, or even the author's own
institutional repository \emph{if} the institution has a
posting mandate.  The last is not a typo: if your institution
mandates posting the accepted author manuscript in its repository, then Elsevier
stipulates that you may not, although they
permit such posting when there is no mandate!

Fortunately, since hearing complaints from the boycotters about
their posting policy, Elsevier has introduced an exception to
allow posting to the arXiv.  However, that is not enough. There
are other non-commercial subject repositories that are
important to segments of the community (Optimization Online,
the Cryptology ePrint Archive, etc.), and more will undoubtedly
be created in the future. Elsevier should allow authors to post accepted
manuscripts to any such repository, as well as to
university repositories, regardless of whether there is a posting mandate.
Furthermore, this right should be guaranteed
by the publishing agreement, not just by a posting policy that
is subject to change at any time.

\subs{Ethics and Peer Review.} Another source of frustration
with Elsevier is their history of lapses in peer review and
ethics. The case of the journal \emph{Chaos, Solitons \& Fractals} (CS\&F) has
become widely known. This journal published 273 papers by its
own editor in chief over eighteen years, 57 of them in a single
year. Suspicions that these papers were not subject to peer
review are corroborated by the editor's declaration that
``senior people are above
this childish, vain practice of peer review.''%
\footnote{C.~Whyte, El Naschie questions journalist
    in Nature libel trial, updated November 16, 2011,
    \url{http://www.newscientist.com/article/dn21169}.}
Elsevier owes the community an explanation for this and other fiascos.
Was there no oversight in place?  Have changes been made so this will
not happen again? What about the other papers in CS\&F?  Are there
records of peer review?  Will any papers that were not peer reviewed
be retracted, or otherwise flagged?  The current situation leaves the
literature in a bad state and compromises the position of authors who
submitted manuscripts for peer review in good faith.  If Elsevier wants
to place this issue behind them, they need to deal with it thoroughly,
forthrightly, and transparently.

In another notorious case, for five years Elsevier ``published
a series of sponsored article compilation publications, on
behalf of pharmaceutical clients, that were made to look like
journals and lacked the proper
disclosures.''\footnote{Statement from Michael
    Hansen, CEO of
    Elsevier's
    Health Sciences Division, regarding Australia based sponsored
    journal practices between 2000 and 2005, May 7, 2009,
    \url{http://www.elsevier.com/wps/find/authored_newsitem.cws_home/companynews05_01203}.}

There are other incidents in which peer review has failed at Elsevier journals,
sometimes in spectacular fashion.%
\footnote{\url{http://umn.edu/~arnold/reasons.html}.} For many
of us, these call into question Elsevier's ability to meet the
standards of quality and ethics we require if we are to
collaborate with them.

\subs{Initial Responses to the Boycott.} On February 27,
Elsevier publicly withdrew its support for the Research Works
Act, which would have prohibited open access mandates for
government-funded research.  The bill was declared dead by
its sponsors in Congress on the very same day. This
victory confirmed the boycott's success in delivering a message
where we were never able to get through before.

Further confirmation came that day in an open letter  from Elsevier
senior vice presidents David Clark and Laura Hassink
to the mathematics community.\footref{fn} %
Besides reporting the about-face on the Research Works Act,
they announced the target price for ``core mathematics titles''
that we discussed above.  They also stated, correctly, that it
would be necessary to address concerns about ``large discounted
agreements'' (bundling) and said that this will come.

Finally, Clark and Hassink announced that free access has been
granted to the archives of fourteen core mathematics journals for the
years from 1995 through four years before the present day.
Access to back issues is indeed critical, and we strongly
believe that all research papers should be made freely
available long before copyright expires. The shorter the delay
the better, of course, but we consider four years a defensible
choice, compatible with the subscription model for journal publishing. The AMS's
experience with a five-year window shows that such a move is
financially viable. We hope that Elsevier's announcement is
just the first step and that expansion to the full set of
mathematical journals and the period before 1995 will be
announced soon.%
\footnote{All three journals discussed
here began publishing in the 1960s.  The issues before 1995 are
currently available from Elsevier online but remain behind
their paywall.}
We also
hope that this is not just a temporary measure. A binding
commitment not to revoke access in the future would
be reassuring on that point.

\subs{Moving Forward.}  While the mathematical literature
itself is a treasure, the current system of scholarly
publishing is badly broken. Elsevier is the largest and, in our
view, the most egregious example of what is wrong. We hope many
readers will agree with us that by choosing to withdraw our
cooperation from Elsevier, we are sending a valuable message to
them and to the scholarly publishing industry more broadly.
Please consider joining the movement at
\url{http://thecostofknowledge.com}.

What is our vision for the future?  The mathematical community
needs a period of experimentation and healthy competition, in
which a variety of publishing models can flourish and develop.
Possibilities include various approaches to open access
publishing,%
\footnote{For example, based on publication charges or on sponsoring
consortia such as SCOAP$^3$ (\url{http://scoap3.org}).}
refereed journals tightly integrated with the arXiv or similar
servers, increased reliance on non-profit publishers, hybrid
models in which community-owned journals subcontract their
operations to commercial publishers, commercially owned
journals with reasonable prices and policies, etc.  It is too early
to predict the mix of models that will emerge as the most
successful.  However, any publisher that wants to be part of this mix
must convince the community that	
they oversee peer review with integrity, that they aid
dissemination rather than hinder it, and that they work to make high-quality
mathematical literature widely available at a reasonable price.

Let's work together to foster good practices and build better
models.
The future of mathematics publishing is in our hands.
\end{document}